\documentclass[12pt]{article}

\usepackage{amssymb,latexsym,amsthm,amsfonts}

\usepackage{graphicx}

\newtheorem{theorem}{Theorem}
\newtheorem{lemma}{Lemma}

\newtheorem{cor}{Corollary}

\newtheorem*{thm}{Theorem}

\begin{document}
    
\title {Multi-Model Cantor Sets}    
\author{Elizabeth Cockerill}
\date{January 21, 2002}
\maketitle

\begin{abstract}	
	In this paper we define a new class of metric spaces, called
	multi-model Cantor sets.  We compute the Hausdorff dimension
	and show that the Hausdorff measure of a multi-model Cantor
	set is finite and non-zero.  We then show that a bilipschitz
	map from one multi-model Cantor set to another has constant
	Radon-Nikodym derivative on some clopen.  We use this to
	obtain an invariant up to bilipschitz homeomorphism.
\end{abstract}

\section*{Introduction}

\noindent A multi-model Cantor set is a metric space which
has the following property.  There is a partition of $C$ into finitely
many clopens $A_{1}, A_{2}, \ldots, A_{n}$ (called models) so that
given any point $x$ in $C$ and $\epsilon >0$ there is a neighborhood
$U$ of $x$ with $diam(U) < \epsilon$ and a metric similarity mapping
$U$ onto one of the models.

Such a Cantor set is described by a map $\tau : C \rightarrow C$ 
which is a piecewise metric similarity (an expanding map).  One may 
consider the Cantor set as determined by the dynamical system $\tau$.  
The middle third Cantor set is an example of a multi-model Cantor set 
where one model suffices.

The middle third Cantor set is self-similar.  At every scale it
contains identical copies of itself.  In a sense a multi-model Cantor
set has finitely many local pictures which are replicated, at
different scales, everywhere.  This behaviour is similar to certain 
fractals, such as Julia Sets, which have a compact family of local 
geometries, (rather than a finite family).

We associate to $C=C(\tau)$ and a constant $d>0$, an $n \times n$ matrix 
$M_{d}$ with non-negative entries.  Roughly speaking the $i,j$ entry 
of $M_{d}$ is the sum of the $d$-powers of the inverse scale 
factors of the similarities of the clopens (later called clones) of 
type $i$ contained in model $j$.  This matrix determines the Hausdorff 
measure and dimension as follows.  We use ${\cal H}_{\delta}$ to 
denote the $d$-dimensional Hausdorff measure.\\

\noindent {\bf Theorem 1.}
{\it Suppose $d>0$ and $C$ is a multi-model Cantor set with matrix
$M_{d}$ and $\lambda_{d}$ is the Frobenius eigenvalue for each
$M_{d}$.  Let $d$ be chosen such that $\lambda_{d}=1$.  Then the
Hausdorff measure of $C$ is finite and non-zero in this dimension. 
Therefore $d$ is the Hausdorff dimension of $C$.  Also, let $\vec{v}$ be the
Frobenius eigenvector of $M_{d}$ such that $\sum \limits_{i=1}^{n}
v_{i} = {\cal H}_{d}(C)$ then $v_{i}= {\cal H}_{d}(A_{i})$.}\\

This result is given at the end of section \ref{Hausdorff}.  Then in
section \ref{Linearity} we investigate bilipschitz maps between
multi-model Cantor sets.  We prove that every bilipschitz map from one
multi-model Cantor set to another is measure linear (has constant
Radon-Nikodym derivative) on some clopen.  This is a generalization of
the results of Cooper \cite{C}, Pignataro \cite{P} and Vu'o'ng \cite{V},
\cite{H} to this wider class of Cantor set.  We use this in corollary 
\ref{10.7} to provide an invariant up to bilipschitz homeomorphism of 
such Cantor sets.  These results are from the author's Master's 
thesis \cite{L}.

\section{Definitions}
\label{Measures}

We are concerned, in this paper, with the study of Cantor sets with
particular metrics.  In the following sections we shall show some
results about the Hausdorff measure of these Cantor sets.  In order to
make these results clearer to the reader we include the definitions of
similarity map, $K$-bilipschitz map, and diameter and state two
elementary results concerning Hausdorff dimension.\\


    \noindent {\bf Definition.} A surjective map $f:X \rightarrow Y$
    between metric spaces is a {\em similarity map} if there is a
    constant $K>0$ such that for all $x_{1}, x_{2}$ in $X$
    \[d_{Y}(f(x_{1}), f(x_{2})) = Kd_{X}(x_{1},x_{2}).\]
    

    \noindent {\bf Definition.} A map $f:X \rightarrow Y$ between
    metric spaces is {\em $K$-bilipschitz} if there is a constant
    $K>0$ such that for all $x_{1}, x_{2}$ in $X$
    \[\frac{1}{K}d_{X}(x_{1},x_{2}) \leq d_{Y}(f(x_{1}), f(x_{2}))
    \leq Kd_{X}(x_{1},x_{2}).\]
    
    
\noindent {\bf Definition.} For a set $U$ contained in a metric
    space $X$, the {\em diameter} of $U$ is given by $$diam(U) =
    \sup\limits_{x,y \in \, U}d_{X}(x,y).$$

\begin{lemma} \label{A}
    Suppose $X$ is a metric space and $d > 0$ and the 
    $d$-dimensional Hausdorff measure of $X$ satisfies $0 < {\cal 
    H}_{d}(X) < \infty$.  Then $d$ is the Hausdorff dimension of $X$.
\end{lemma}

\begin{proof}
    This follows directly from the definition.
\end{proof}

\begin{lemma}
    \label{multiply}
    Let $X$, $Y$ be metric spaces.  Let $f:X \rightarrow Y$ be a
    surjective similarity map with scale factor $K>0$, then $${\cal
    H}_{\delta}(Y) =
    K^{\delta}{\cal H}_{\delta}(X).$$
\end{lemma}

\begin{proof}
    This follows directly from the definition.
\end{proof}

\section{The Middle Third Cantor Set}
\label{middle third}

Recall that a topological space $C$ is a {\em Cantor set} if it is
compact, totally disconnected, perfect and non-empty.  It is a fact
that any two Cantor sets are homeomorphic.  The situation becomes more
interesting when one puts a metric on a Cantor set and studies
bilipschitz homeomorphisms, as shall be seen in section
\ref{Linearity}.

In this paper we are interested in considering Cantor sets which
arise from dynamical systems.  We illustrate this by exhibiting a dynamical
system which generates the middle third Cantor set.  This dynamical
system is determined by a {\em clone structure} which we shall now
describe.

Let $I=[0,1]$, $A_{1}=\left[0,\frac{1}{3}\right]$ and
$A_{2}=\left[\frac{2}{3},1\right]$ and define the {\em clone maps}
$\tau_{i}:A_{i} \rightarrow I$ by $$\tau_{1}=3x \quad \mbox{and} \quad
\tau_{2}=3x-2.$$ Then this {\em clone structure} is either denoted as
$(A_{1}, \tau_{1}, A_{2}, \tau_{2})$ or, because the domain $A_{i}$ is
implicit in the definition of $\tau_{i}$, it can also be denoted
$(\tau_{1}, \tau_{2})$.

To construct the middle third Cantor set from the dynamical system we
define the map $\tau = \tau_{1} \cup \tau_{2}:A_{1} \cup A_{2}
\rightarrow I$ by $$ \tau(x) = \left\{ \begin{array}{ll} 3x & \mbox{
if $x \in \,  A_{1}$}; \\ 3x-2 & \mbox{ if $x \in \,  A_{2}$}.  \end{array}
\right.$$ Observe that $\tau$ maps each of the $A_{i}$
linearly onto $I$.  Then $$\tau^{-1}(I) = A_{1} \cup A_{2} = C_{1}.$$
In general $$\tau^{-n}(I)= C_{n}$$ is a disjoint union of $2^{n}$
intervals, each of length $3^{-n}$.  The middle third Cantor set is
$$C(\tau) = \bigcap\limits_{n=0}^{\infty}\tau^{-n}(I).$$

\begin{figure}[htbp]
\begin{center}
\includegraphics[width=105mm,height=20mm]{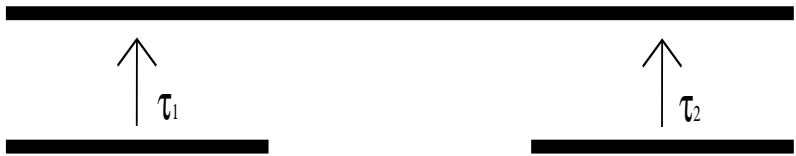}
\end{center}
\caption{\label{middlethird3}Middle third Cantor Set with clone structure 
$(\tau_{1}, \tau_{2})$.}
\end{figure}

\begin{figure}[htbp]
\begin{center}
\includegraphics[width=105mm,height=20mm]{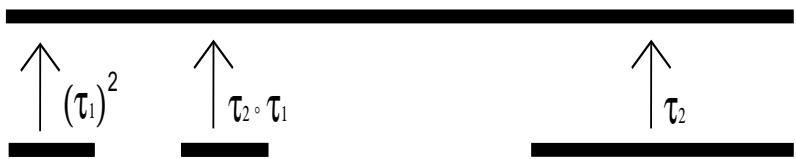}
\end{center}
\caption{\label{middlethird4}Middle third Cantor Set with clone structure 
$\left(\tau_{1} \circ \tau_{1}, \tau_{2} \circ \tau_{1}, \tau_{2}\right)$.}
\end{figure}

We must note though that the clone structure is not unique.  At each
level any clone can be replaced by a finite collection of clones which
partition it.  Therefore we usually refer to a Cantor set $C$, rather
than $C(\tau)$.  In figures \ref{middlethird3} and \ref{middlethird4}
we show two different clone structures for the middle third Cantor set.

Having chosen a clone structure we shall call a clopen a {\em clone}
if it is the domain of a composition of clone maps in the clone
structure.  The empty composition is the identity map so the entire
Cantor set $C$ is also a clone.  The {\em level} of the Cantor set is
defined to be 0.  The {\em level} of any other clone is defined to be
$n$ if the smallest clone which properly contains it has level $n-1$. 
Hence the middle third Cantor set has two level-1 clones: $A_{1}=C \,
\cap \, \left[0, \frac{1}{3}\right]$ and $A_{2} = C \, \cap \,
[\frac{2}{3},1].$ Note that we will abuse notation and use $A_{1}$ to
represent both the interval $[0, 1/3]$ and also the subset of the
Cantor set contained in $[0, 1/3]$.  A clone $B$ is said to be at {\em
relative level 1} to a clone $D$ if given the level of clone $D$ as
$n$ then the level of clone $B$ is $n+1$.

The middle third Cantor set is generated by 2 similarity maps, each with a
similarity factor of 3.  Cantor sets can be generated by two or more maps
with various similarity factors.  They can also be generated in higher
dimensions in Euclidean space \cite{V}, \cite{H} or in a general metric space
\cite{C}.  Below is a pictorial representation of a Cantor set in
$\mathbb{R}^{2}$.  This Cantor set is given by a dynamical
system; a clone structure.

\begin{figure}[htbp]
\begin{center}
\includegraphics[height=45mm]{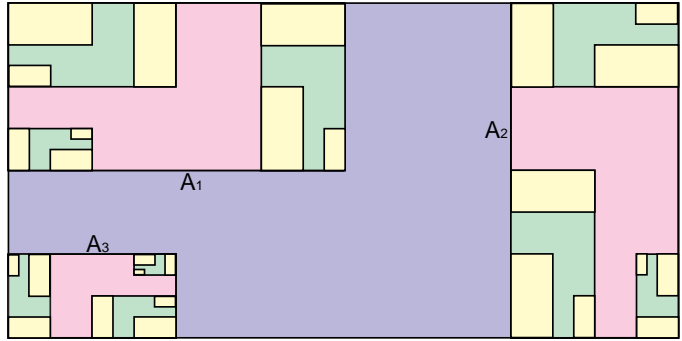}
\end{center}
\caption{\label{2dcantor5}Cantor set in $\mathbb{R}^{2}$ with clone 
structure $(\tau_{1}, \tau_{2}, \tau_{3})$.}
\end{figure}

Figure \ref{2dcantor5} is a Cantor set in $\mathbb{R}^{2}$.  It has 3
level-1 clones, $A_{1}, A_{2}, A_{3}$.  Below we describe how each
level-1 clone is mapped onto the largest rectangle.  Each clone
contains a scaled down copy of the entire Cantor set.  A clone
structure on this Cantor set can be described by three similarity
maps.  If the origin is considered to be the top left vertex of the
rectangle $A$, then the map $\tau_{1}:A_{1} \rightarrow A$ is a
stretch by a factor of 2, $\tau_{2}:A_{2} \rightarrow A$ is a rotation
by $\frac{\pi}{2}$ radians clockwise, a stretch by a factor of 2 and a
translation and $\tau_{3}:A_{3} \rightarrow A$ is a rotation by $\pi$
radians, a stretch by a factor of 4 and a translation.

\section{Multi-Model Cantor Sets}
\label{multi-model}


In section \ref{middle third} we discussed clone Cantor sets in metric
spaces arising from dynamical systems, where each clone can be mapped
onto the whole Cantor set by using a similarity clone map.  In this
section we shall be looking at a new class of examples, called
{\em multi-model Cantor sets} which also arise from dynamical systems.

In Vu'o'ng's thesis \cite{V} and paper \cite{H}, he studied clone
Cantor sets which can be determined by collections $\{\tau_{1},
\tau_{2}, \ldots, \tau_{q}\} $ where the $A_{i}$ are level-1 clones of
$C$ and the $\tau_{i}:A_{i} \rightarrow C$ are level-1 clone maps of
$C$.  Vu'o'ng was considering Cantor sets where each clone is mapped
onto the unit cube in $\mathbb{R}^{n}$ by a clone map.  In other words
he was considering Cantor sets which look the same everywhere.  In a
sense a multi-model Cantor set has finitely many local pictures which
are replicated, at different scales, everywhere.\\

\noindent {\bf Definition.} Let $(C,d)$ be a metric space.  Let
$A_{1}, A_{2}, \ldots, A_{n}$ be a partition of $C$ by clopens.  The
$A_{j}$ are called {\em models}.  Let $E_{1}, E_{2}, \ldots, E_{m}$ be
another partition of $C$ by clopens such that each $E_{i}$ is
contained in some $A_{j}$.  The $E_{i}$ are called {\em level-1
clones}.  We require that each $A_{j}$ contains at least 2 of the
$E_{i}$ hence $m \geq 2n$.  For each $E_{i}$ there is a surjective
similarity map $\tau_{i}:E_{i} \rightarrow A_{j}$ for some $j$.  We
say that $(\tau_{1}, E_{1},\tau_{2}, E_{2}, \ldots, \tau_{m}, E_{m})$
is a {\em clone structure} of a {\em multi-model Cantor set}, $C$.

As before we define the piecewise metric similarity $$\tau = \bigcup
\limits_{i=1}^{m} \tau_{i} : \bigcup \limits_{i=1}^{m} E_{i}
\rightarrow \bigcup \limits_{j=1}^{n} A_{j}.$$\\
We may construct examples in Euclidean space by using similarities of
$\mathbb{R}^{n}$.  For example if $A$ is a hypercube in
$\mathbb{R}^{n}$ we obtain a nested set of clopens $\{ C_{n} \}$ where
$$C_{n}= \tau^{-n}(A).$$
Then the multi-model Cantor set is the intersection $$C(\tau)=\bigcap
\limits_{n=0}^{\infty} \tau^{-n}(A).$$

Each model and each clone in a multi-model Cantor set has a {\em type}.  By
definition each model $A_{j}$ is of type $j$.  Each clone $E_{i}$ is
defined to be the same type as the model it is mapped to by
$\tau_{i}$.  Clones at higher levels are defined to be of the same
type as the model they are mapped to by compositions of the similarity
maps, $\tau_{i}$.

Each $\tau_{i}$ is a similarity map from a level-1 clone to a model. 
There is a lot of information we need to convey with the label
$\tau_{i}$, such as which model the clone $E_{i}$ is mapped to and
which it is contained in.  For the sake of clarity we use only one
subscript throughout most of this paper.  In the next two paragraphs
we include an explanation of the rules for compositions of these
similarity maps.  To do this we introduce three indices for each map
which we shall subsequently omit in later sections.

We define $\tau_{i}=\tau_{\alpha, \beta, r}:E_{i} \to A_{\beta}$ where
$1 \leq \alpha \leq n$, $1 \leq \beta \leq n$ and $1 \leq r \leq
R_{\alpha,\beta}$.  Here $\alpha$ is the type of model the level-1
clone $E_{i}$ is contained in, $\beta$ is the type of the level-1
clone $E_{i}$, and $R_{\alpha,\beta}$ is the number of level-1 clones
of type $\beta$, contained in model $\alpha$.

\begin{figure}[htbp]
\begin{center}
\includegraphics[width=105mm,height=30mm]{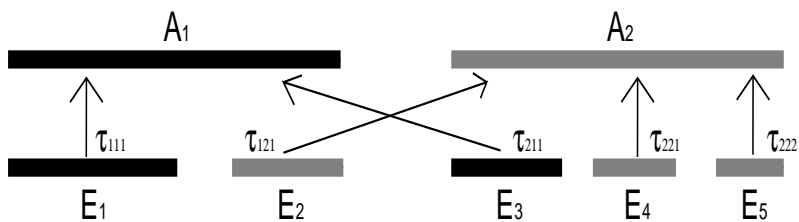}
\end{center}
\caption{\label{cantor11}Models, level-1 clones and similarity maps. 
(Type 1 is dark, type 2 is light.)}
\end{figure}

The compositions of these similarity maps are called {\em clone-model
maps}.  They carry the restriction that if $\tau_{\alpha, \beta, r}$
is composed with $\tau_{\gamma, \delta, l}$ to form $\tau_{\gamma,
\delta, l}\circ \tau_{\alpha, \beta, r}$ then $\alpha$ must be equal
to $\delta$.  For example $\tau_{1,{\bf 3},1} \circ \tau_{{\bf
1},2,2}$ is not a possible composition, but $\tau_{2,{\bf 3},1}\circ
\tau_{{\bf 3},2,2}$ is.  Also, the index in the $\beta$ position of
the first map and the index in the $\alpha$ position on the last map
must be the same.  For example, $\tau_{{\bf 2},3,1} \circ \tau_{3,1,1}
\circ \tau_{1,{\bf 3},1}$ is not defined but $\tau_{{\bf 3},3,1} \circ
\tau_{3,1,1} \circ \tau_{1,{\bf 3},1}$ is defined.

The domain of each clone-model map is by definition a {\em clone}.  A
{\em level-k} clone is the domain of a composition of $k$ of
the clone-model maps.

Below are three pictorial representations of a single multi-model
Cantor set in $\mathbb{R}^{2}$.  This multi-model Cantor set is given
by the dynamical system $\tau = (\tau_{1}, E_{1}, \tau_{2}, E_{2},
\ldots, \tau_{8}, E_{8})$.

\begin{figure}[htbp]
\begin{center}
\includegraphics[height=55mm]{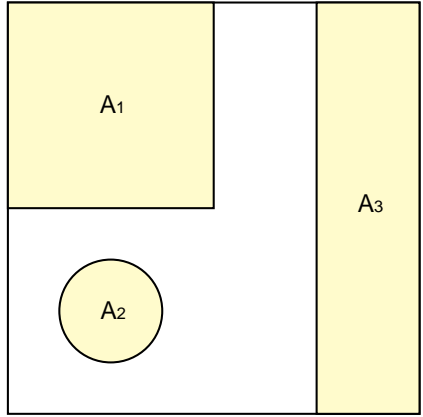}
\end{center}
\caption{\label{2dmulti6}The models of a multi-model Cantor set (shown
fully in figure \ref{2dmulti4}).}
\end{figure}

\begin{figure}[htbp]
\begin{center}
\includegraphics[height=55mm]{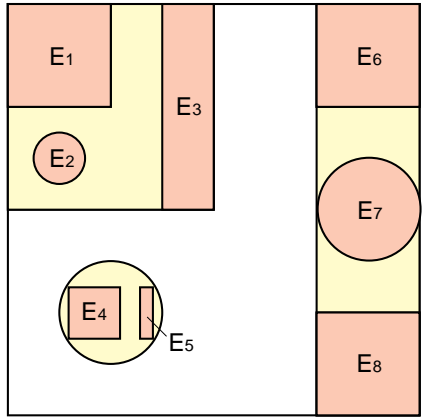}
\end{center}
\caption{\label{2dmulti5}The models and level-1 clones of figure 
\ref{2dmulti4}.}
\end{figure}
    
\begin{figure}[htbp]
\begin{center}
\includegraphics[height=55mm]{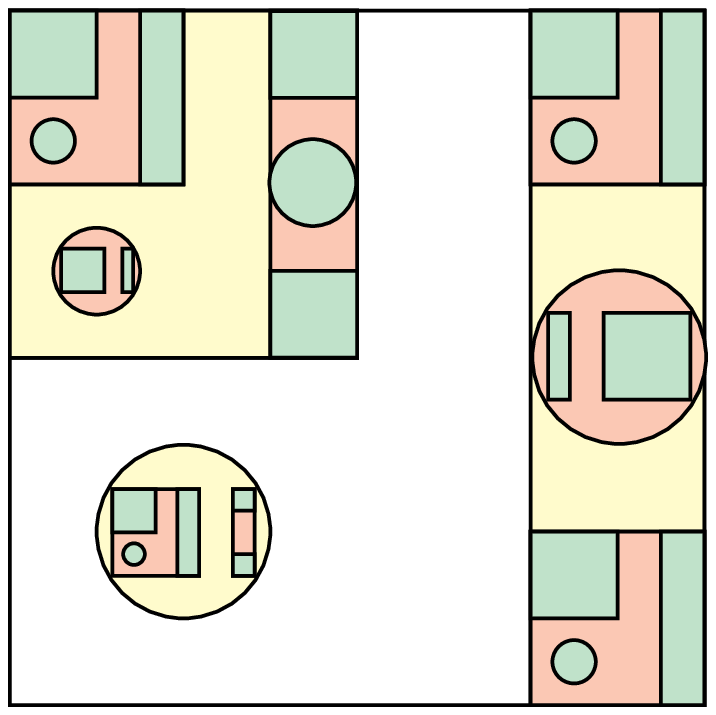}
\end{center}
\caption{\label{2dmulti4}A multi-model Cantor set in $\mathbb{R}^{2}$ 
with clone structure $(\tau_{1}, E_{1}, \tau_{2}, E_{2}, \ldots, 
\tau_{8}, E_{8})$.}
\end{figure}

Figure \ref{2dmulti4} is a multi-model Cantor set in $\mathbb{R}^{2}$. 
The models, seen in figure \ref{2dmulti6}, are $A_{1}$ - a square,
$A_{2}$ - a circle and $A_{3}$ - a rectangle.  The level-1 clones
$E_{1}, E_{2}, \ldots, E_{8}$ are contained in the models.  These can
be seen more clearly in figure \ref{2dmulti5}.

We can see that the similarity maps from the level-1 clones onto the
models are quite simple.  If we assume the top left vertex to be the
origin then clone $E_{1}$ is mapped onto $A_{1}$ by a stretch factor
of 2.  Clone $E_{2}$ maps onto $A_{2}$ by a stretch factor of 2 and a
translation.  Clone $E_{3}$ maps onto model $A_{3}$ by a stretch
factor of 2 and a translation.  Clone $E_{4}$ maps onto $A_{1}$ by a
stretch factor of 4 and a translation.  Clone $E_{5}$ maps onto
$A_{3}$ by a stretch of factor 4 and a translation.  Clone $E_{6}$
maps onto $A_{1}$ by a stretch of factor 2 and a translation.  Clone
$E_{7}$ maps onto $A_{2}$ by a rotation by $\pi$ and a translation and
$E_{8}$ maps onto $A_{1}$ by a stretch of factor 2 and a
translation.\\

    \noindent {\bf Definition.} Given a multi-model Cantor set
    $C=C(\tau)$ with $n$ models, we define for each $d>0$ an $n \times
    n$ matrix $M_{d}$ as follows.  The $i,j$ entry is the sum of the
    $d$-powers of the inverse scale factors of the clone maps of the
    level-1 clones of type $i$, contained in model $j$.  Recall,
    $R_{ij}$ is the number of type $i$ clones contained in model $j$. 
    So each entry is of the form $(M_{d})_{i,j}= \sum\limits
    ^{R_{ij}}_{r=1} \left((a_{ij})_{r} \right)^{d}$ where $0 <
    (a_{ij})_{r} <1$ is the inverse of the scale factor of the
    revelant clone map.
    
\begin{figure}[htbp]
\begin{center}
\includegraphics[width=105mm,height=35mm]{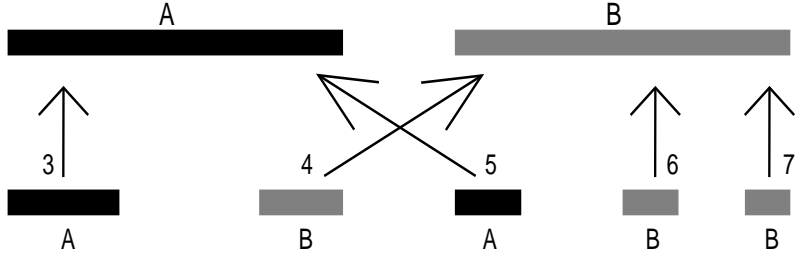}
\end{center}
\caption{\label{matrix}A multi-model Cantor set and its 
corresponding matrix.}

\[ M_{d}={\large \left(\begin{array}{cc}\left (\frac{1}{3}
\right)^{d}& \left (\frac{1}{5} \right)^{d} \\
    \left (\frac{1}{4} \right)^{d}& \left
(\frac{1}{6} \right)^{d} +  \left
(\frac{1}{7} \right)^{d}\end{array} \right).} \]\\
\end{figure}

Figure \ref{matrix} gives an example of a multi-model Cantor set and 
its corresponding matrix.  The scale factors of the similarity maps 
are given on the diagram.

\noindent {\bf Definition.} Given a multi-model Cantor set $C$ with
$n$ models and a finite collection of disjoint clones ${\cal J} = \{
J_{i,r_{i}} \} $, where $i$ is the type of clone and $1 \leq r_{i}
\leq R_{i}$ where $R_{i}$ is the number of type $i$ clones in {\cal J}. 
Then given $d>0$ we define a vector $\vec{v} = (v_{1}, \ldots
v_{n})^{T}$ by $$v_{i} = \sum \limits_{r_{i}=1}^{R_{i}}
(diam(J_{i,r_{i}}))^{d}.$$ We say $\vec{v} = \vec{v}({\cal J})$ is the
{\em d-quantity} of the collection of clones ${\cal J}$.\\

\noindent {\bf Definition.} Given a finite collection of disjoint
clones ${\cal J}$ and $k \geq 1$ let ${\cal J}^{(k)}$ be the
collection of clones obtained by subdividing, $k$ times, each clone in
${\cal J}$.  (By subdividing once we mean replacing each clone by the
clones contained in it at relative level 1.)

\begin{lemma}
    \label{power}
    Given a multi-model Cantor set $C=C(\tau)$ and given $d>0$ let
    $M=M_{d}$ be the the multi-model Cantor set matrix defined above. 
    Given a finite collection of clones ${\cal J}$ then
    $$\vec{v}\left({\cal J}^{(k)}\right) = M^{k}\, \vec{v}\left({\cal
    J}\right ).$$

\end{lemma}

\begin{proof}
    
    First observe that if ${\cal J} = {\cal J}_{1} \cup {\cal J}_{2}$
    and ${\cal J}_{1} \cap {\cal J}_{2} = \emptyset$ then
    $\vec{v}\left({\cal J}\right) = \vec{v}\left({\cal J}_{1}\right) +
    \vec{v}\left({\cal J}_{2}\right).$ Thus it suffices to prove the
    result when ${\cal J}$ is a single clone.  The case $k=1$ follows
    from the definition of $M_{d}$ (refer to figure \ref{matrix}). 
    The inductive step follows from the observation that
    $${\cal J}^{(k+1)} = \left ({\cal J}^{(k)}\right )^{(1)}$$  
    and applying the case $k=1$ to the collection ${\cal J}^{(k)}.$
    \end{proof}

\begin{cor}
    \label{mtothek}
    Given a multi-model Cantor set $C(\tau)$ with multi-model Cantor 
    set matrix $M_{d}$, then $(M_{d})^{k}$ is the multi-model Cantor 
    set matrix for $C(\tau^{k}).$  (Observe that the Cantor sets 
    $C(\tau)$ and $C(\tau^{k})$ are identical.)
\end{cor}

    \noindent {\bf Definition.} A square matrix $M$ is {\em
    irreducible} if all its entries are non-negative and there exists
    some $k>0$ such that all entries of $M^{k}$ are positive.  (Note
    this is not the same definition of irreducible used in
    representation theory; namely that there are no invariant
    subspaces).

In this paper we study multi-model Cantor sets with
irreducible matrices.  The geometric significance of irreducibility, by
corollary \ref{mtothek}, is that a clone of each type is contained at
some level within each model.  Here are two figures which illustrate
this.  The first is irreducible, the second is not.

\begin{figure}[htbp]
\begin{center}
\includegraphics[width=105mm,height=30mm]{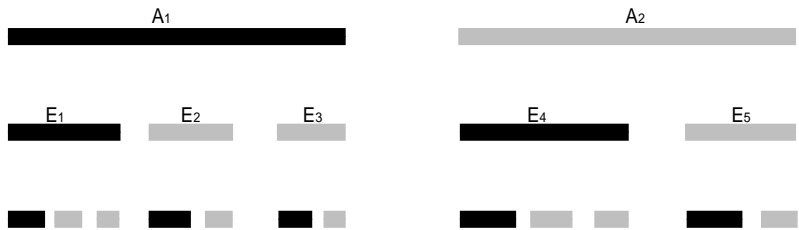}
\end{center}
\caption{\label{irreducible}A multi-model Cantor set with irreducible
matrix.  (Each clone is contained in each model.)}
\end{figure}

\begin{figure}[htb]
\begin{center}
\includegraphics[width=105mm,height=30mm]{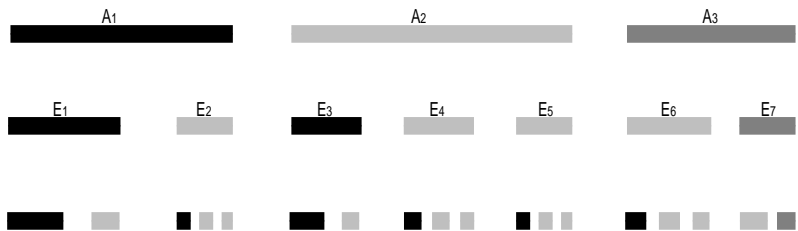}
\end{center}
\caption{\label{reducible}A multi-model Cantor set with non-irreducible
matrix.  (Type 3 clones are not contained in models $A_{1}$ or $A_{2}$.}
\end{figure}

Figure \ref{irreducible} shows a multi-model Cantor set with an
irreducible matrix.  There are no zeros in the corresponding matrix,
$(M_{d})^{k}$, for any $k>0$.

In figure \ref{reducible} we see a multi-model Cantor set which does
not have an irreducible matrix.  The corresponding matrix
$(M_{d})^{k}$ has zero entries in the $3rd$ row of
both the $1st$ and $2nd$ columns, for every $k>0$.

\section{The Hausdorff Measure Of Multi-Model Cantor Sets}
\label{Hausdorff}

In this section we prove theorem \ref{finite measure}.

\begin{thm}[{\bf Frobenius}]
    \label{Frob}
    An $n \times n$ irreducible matrix $M$ always has a real
    positive eigenvalue $\lambda$ that is a simple root of the
    characteristic equation and is larger than the moduli of all the other
    eigenvalues.  To this `maximal' eigenvalue
    $\lambda$ there corresponds an eigenvector with strictly positive
    coordinates.
\end{thm}

    \noindent {\bf Definition.} The {\em Frobenius eigenvalue} of an
    irreducible matrix is the eigenvalue of largest modulus.  From the
    Frobenius theorem we see that it is positive, real and has
    multiplicity 1.  A {\em Frobenius eigenvector} is a corresponding
    positive, real eigenvector.
    
    The first step in the proof of theorem \ref{finite measure} is to 
    show that the Frobenius eigenvalue $\lambda_{d}$, of the 
    multi-model Cantor set matrix $M_{d}$, is a strictly decreasing 
    continuous function of $d$.  From this we shall see that there is 
    a unique $d>0$ such that $\lambda_{d}=1$, and that this is the 
    Hausdorff dimension of the Cantor set.  We state some easily 
    proved results on matrices without proof.\\

\noindent {\bf Definition.}  The {\em $L_{1}$-norm} on 
$\mathbb{R}^{n}$ is $\| \, \, \|_{1}:\mathbb{R}^{n} \rightarrow 
\mathbb{R}$, given by $$\| (v_{1}, \ldots, v_{n}) \|_{1} = \sum 
\limits_{i=1}^{n} |v_{i}|.$$

\begin{lemma}
    \label{prop 1}
    Let $M$ be an $n \times n$ irreducible matrix with Frobenius
    eigenvalue $\lambda$.  Then there exists $Q>0$ such
    that for all $k$ and for all vectors $\vec{v}$,
    $$\|M^{k}\vec{v}\|_{1} \leq Q\lambda^{k}\|\vec{v}\|_{1}.$$ Furthermore
    when $\lambda=1$ there is a strictly positive matrix
    $M^{\infty}$ such that $M^{k}$ converges to $M^{\infty}$ as $k$
    tends to $\infty$.
\end{lemma}

    \noindent {\bf Definition.} We put a partial ordering on
    real matrices.  We say matrix $A$ is {\em less than} matrix $B$ if
    each entry in $A$ is less than the corresponding entry in $B$.  So
    $A<B$ if $A_{ij} < B_{ij}$ for all $i$ and $j$.
    
\begin{lemma}
    \label{eigen}
    If $\lambda$ is an eigenvalue of $M$ with 
    eigenvector $\vec{v}$ then $\lambda^{k}$ is an eigenvalue of 
    $M^{k}$ with eigenvector $\vec{v}$.
\end{lemma}

Note, if $\lambda$ is the Frobenius eigenvalue then it has largest
modulus of all eigenvalues for $M$ so $\lambda^{k}$ has largest
modulus of all eigenvalues for $M^{k}$.  Therefore $\lambda^{k}$ is
the Frobenius eigenvalue for $M^{k}$.

\begin{lemma} \label{bigger}
    Suppose $0 < M_{1}< M_{2}$, and $\lambda_{i}$ is the Frobenius
    eigenvalue for $M_{i}$, then $\lambda_{1}<\lambda_{2}$.
\end{lemma}

\begin{cor}
    \label{decreases}
    Suppose $M_{d}$ is an $n \times n$ strictly positive matrix such
    that each entry is a strictly monotonically decreasing continuous function
    of $d$.  Let $\lambda_{d}$ be the Frobenius eigenvalue of $M_{d}$
    for each $d$.  Then as $d$ increases, $\lambda_{d}$ decreases.
\end{cor}

\begin{proof}
    Let $0< d'<d$.  Each entry of $M_{d}$ is a strictly monotonically
    decreasing function of $d$ so $M_{d'} > M_{d} > 0$.  Then by lemma
    \ref{bigger} we have $\lambda_{d'}>\lambda_{d}$.
\end{proof}

\begin{lemma}
    Suppose $M_{d}$ is an $n \times n$ strictly positive matrix where
    each entry is a strictly monotonically decreasing continuous function of
    $d$, and $\lambda_{d}$ is the Frobenius eigenvalue of $M_{d}$ for
    each $d$.  Then $f:[0,\infty) \rightarrow \mathbb{R}^{+}$ such
    that $f(d)=\lambda_{d}$ is continuous.
\end{lemma}

\begin{proof}
    Let $g_{d}(x)$ be the characteristic polynomial of $M_{d}$.  The
    coefficients of $g_{d}(x)$ are polynomial functions of the entries
    of $M_{d}$.  Since polynomials are continuous and the entries of
    $M_{d}$ are continuous functions of $d$, it follows that the
    coefficients of $g_{d}(x)$ vary continuously with $d$.
    
    Frobenius gives us $\lambda_{d}$ which is a real, positive root of
    $g_{d}(x)$ with multiplicity 1.  Without loss of generality we can
    assume $g_{d}(\lambda_{d}+\epsilon)>0$ and $g_{d}
    (\lambda_{d}-\epsilon) <0$ for some small $\epsilon > 0$.  Now
    consider $g_{d'}(x)$ where $d'=d+\delta$ for some small $\delta > 0$,
    then $g_{d'}(\lambda_{d}+\epsilon)>0$ and
    $g_{d'}(\lambda_{d}-\epsilon)<0$ by continuity of coefficients of
    the characteristic polynomial.  Then by the intermediate value
    theorem $g_{d'}(x_{d'})=0$ for some $\lambda_{d}-\epsilon < x_{d'} <
    \lambda_{d}+\epsilon$.
    
    We have shown that $g_{d'}(x)$ has a root $x_{d'}$ in the interval
    $(\lambda_{d}- \epsilon, \lambda_{d} +\epsilon)$.  We now complete
    the proof that $f(d) = \lambda_{d}$ is continuous.  We argue by
    contradiction.  If as $d' \to d$ we have $\lambda_{d'} \not\to
    \lambda_{d}$ then, since $\lambda_{d}$ is monotonic in $d$, there
    is an $\eta > 0$ and a sequence $\{d_{n}\}$ such that $d_{n} \to
    d$ with $| \lambda_{d} - \lambda_{d_{n}}| > \eta.$ Since
    $\lambda_{d_{n}}$ is the largest root we know that
    $\lambda_{d_{n}} \geq x_{d_{n}}$ thus $\lambda_{d_{n}} >
    \lambda_{d} + \eta$.  For $n$ large we have $d-1 \leq d_{n} \leq
    d+1$ so by corollary \ref{decreases} we have $\lambda_{d+1} \leq
    \lambda_{d_{n}} \leq \lambda_{d-1}$.  Therefore the
    $\lambda_{d_{n}}$ are contained in a compact interval and so there
    exists a convergent subsequence.  Re-label the convergent
    subsequence, then $\lim\limits_{n \to \infty} \lambda_{d_{n}} =
    \lambda' \geq \lambda_{d} + \eta.$
    
\noindent {\bf Claim.}
    The limit $\lambda'$ is a root of $g_{d}(x)$.

\noindent {\em Proof.}
    The sequence $d_{n} \to d$ and so $g_{d_{n}}(x) \to g_{d}(x)$
    because the coefficients of $g_{d}(x)$ are continuous functions of
    $d$.  The characteristic equation $g_{d_{n}}(x)$ has a root
    $\lambda_{d_{n}}$ for all n so $$g_{d}(\lambda') = \lim \limits_{n
    \to \infty} g_{d_{n}}(\lambda_{d_{n}}) = \lim\limits_{n
    \to \infty}(0) = 0.$$
    
    This proves the claim.  The claim contradicts $\lambda_{d}$ being
    the largest root of $g_{d}(x)$.  The contradiction implies that
    $x_{d'}$ is the largest root of $g_{d'}(x)$ and so $x_{d'} =
    \lambda_{d'}$.  Therefore $\lambda_{d}$ varies continuously with
    $d$ and hence $f(d)=\lambda_{d}$ is continuous.
\end{proof}

In the following lemmas we shall suppose that the irreducible matrix
$M_{d}$ is strictly positive.  This can be justified as follows. 
Given $\tau$ and an integer $k >0$ the multi-model Cantor sets
$C(\tau)$ and $C(\tau^{k})$ are identical.  However their clone
structures differ.  The corresponding multi-model Cantor set matrices
are respectively $M_{d}$ and $(M_{d})^{k}$ by corollary
\ref{mtothek}.  The hypothesis that $M_{d}$ is irreducible ensures
that we may choose $k$ sufficiently large so that $(M_{d})^{k}$ is
strictly positive.  Therefore replacing $\tau$ by $\tau^{k}$ justifies
the supposition.

\begin{lemma}
Suppose $M_{d}$ is a multi-model Cantor set matrix and 
$\lambda_{d}$ is the Frobenius eigenvalue for each $d$.  Then there is 
a unique $d>0$ such that $\lambda_{d}=1$.
\end{lemma}

\begin{proof}
    Recall that the inverse scale factors of clone maps satisfy
    $0<(a_{ij})_{r}<1$.  As $d$ approaches $\infty$, we have that
    $\left(M_{d}\right)_{ij} = \sum\limits^{R_{ij}}_{r=1}
    ((a_{ij})_{r})^{d}$ gets close to $0$, for all $i,j$.  This
    implies $M_{d}$ approaches the zero matrix and so $\lambda_{d}$
    tends to $0$ by lemma \ref{bigger}.

    When $d=0$, the sum $\sum \limits_{r=1}^{R_{ij}}
    ((a_{ij})_{r})^{0} = R_{ij}$ for all $i,j$.  Thus each of the
    entries of the matrix $M_{0}$ are at least $1$.  Let $\vec{v}$ be
    the eigenvector for $\lambda_{0}$ then $$\lambda_{0}\vec{v} =
    M_{0}\vec{v} \geq \left(\begin{array}{cccc}1&1&\cdots&1\\
    1&1&\cdots&1\\
    \vdots&\vdots&\ddots&\vdots\\1&1&\cdots&1\\ \end{array}
    \right)\vec{v} = \left(
    \begin{array}{c}\sum\limits^{n}_{i=1}v_{i}\\
    \sum\limits^{n}_{i=1}v_{i}\\ \vdots \\
    \sum\limits^{n}_{i=1}v_{i}\\ \end{array} \right) .$$

    Note, $\vec{v}$ is the Frobenius eigenvector of $M_{0}$ and hence
    $v_{i}>0$ for all $i$, therefore $\sum\limits^{n}_{i=1}v_{i}>
    v_{j}$, for all $j$.  This implies $\lambda_{0}\vec{v} > \vec{v}$
    and hence $\lambda_{0} > 1$.

    The function $f(d)=\lambda_{d}$ is continuous, so by the
    intermediate value theorem $\lambda_{d} = 1$ for some value of
    $d>0$.  The function $f(d) = \lambda_{d}$ is strictly monotonic
    by lemma \ref{bigger}.  Therefore this $d$ is unique.
\end{proof}

\begin{lemma} \label{B}
    Suppose $M$ is an irreducible non-negative $n \times n$ matrix
    with Frobenius eigenvalue $1$.  Then there is a $Q>0$ such that
    for all $p>0$ and for all vectors $\vec{v}$ in $\mathbb{R}^{n}$ we
    have $$\| M^{p} \vec{v} \|_{1} \leq Q \|\vec{v}
    \|_{1}.$$
\end{lemma}

The following gives the Hausdorff measure of every model in a 
multi-model Cantor set.  A clone is a scaled down copy of a model and 
a clopen is a finite union of clones.  Thus the Hausdorff measure of 
every clopen can be determined from the following result.

\begin{theorem}
    \label{finite measure}
Suppose $d>0$ and $C$ is a multi-model Cantor set with matrix
$M_{d}$ and $\lambda_{d}$ is the Frobenius eigenvalue for each
$M_{d}$.  Let $d$ be chosen such that $\lambda_{d}=1$.  Then the
Hausdorff measure of $C$ is finite and non-zero in this dimension. 
Therefore $d$ is the Hausdorff dimension of $C$.  Also, let $\vec{v}$ be the
Frobenius eigenvector of $M_{d}$ such that $\sum \limits_{i=1}^{n}
v_{i} = {\cal H}_{d}(C)$ then $v_{i}= {\cal H}_{d}(A_{i})$.
\end{theorem}

\begin{proof}
    By lemma \ref{A}, if we show that $0<{\cal H}_{d}(C)<\infty$ 
    then $d$ is the Hausdorff dimension of $C$.
    
    To show ${\cal H}_{d}(C) < \infty$ it suffices to show there is 
    a $K>0$ such that for all $\epsilon>0$ there is an open cover 
    ${\cal U}$ of $C$ such that $diam(U) < \epsilon$ for all $U$ 
    in ${\cal U}$ and $$\sum \limits_{U \in \, {\cal U}}(diam(U))^{d} < 
    K.$$
    
    Given $\epsilon >0$ there is a $k>0$ such that the diameter of
    every level-k clone is less than $\epsilon$.  Let ${\cal U}$ be
    the set of level-k clones.  Then ${\cal U}$ is an open
    $\epsilon$-cover of $C$.  If $M=M_{d}$ and if $\{A_{j}\}_{1\leq j 
    \leq n}$ are the
    level-0 models then by lemma \ref{power} the $d$-quantity of 
    ${\cal U}$ is $$\vec{v}({\cal U}) = \sum \limits_{U \in \,
    {\cal U}}(diam(U))^{d} = \sum \limits_{i,j=1}^{n}
    (M^{k})_{ij}(diam(A_{j}))^{d}.$$ Now $M^{k}$ converges to
    $M^{\infty}$ as $k \to \infty$ by lemma \ref{prop 1}.  Hence
    $$\lim \limits_{k \to \infty} \sum \limits_{i,j=1}^{n}
    (M^{k})_{ij}(diam(A_{j}))^{d} = \sum \limits_{i,j=1}^{n}
    (M^{\infty})_{ij}(diam(A_{j}))^{d} = K'.$$ If we take $K=K'+1$
    then we are done.
    
    To show ${\cal H}_{d}(C) >0$ it suffices to show there is an 
    $\eta > 0$ such that for every open cover ${\cal U}$ of $C$ we 
    have $$\sum \limits_{U \in \,  {\cal U}}(diam(U))^{d} > \eta.$$
    
    Given an open set $U$, let $J_{U}$ be the smallest clone containing
    $U$.  By lemma \ref{10.3} there is a constant $\beta =
    \beta(\tau)$ depending only on the Cantor set and $\tau$ such that
    $diam(U) \geq \beta^{-1}diam(J_{U}).$ Hence $$\sum \limits_{U \in
    \, {\cal U}}(diam(U))^{d} \geq \beta^{-d} \sum \limits_{U \in \,
    {\cal U}}(diam(J_{U}))^{d}.$$ Therefore it suffices to show that
    there is an $\eta'>0$ such that for every cover ${\cal U}$ of $C$
    by clones we have \begin{equation} \label{eta} \sum \limits_{U \in
    \, {\cal U}}(diam(U))^{d} > \eta'.  \end{equation} Then $\eta =
    \beta^{-d}\eta'$ satisfies the above.
    
    Given a cover of $C$ by clones, since $C$ is compact there is a
    finite subcover.  We may choose the subcover so that no clone in
    the subcover is contained in any other clone in the subcover.  If
    two clones intersect then one is contained in the other.  Hence we
    may choose a subcover consisting of finitely many disjoint clones. 
    Therefore it suffices to prove inequality (\ref{eta}) for ${\cal
    U}$ a finite cover of $C$ by disjoint clones.
    
    The clones in ${\cal U}$ may be at different levels.  Let ${\cal
    U}_{k}$ be the cover obtained by subdividing the clones in ${\cal
    U}$ so that they are all at the same level, $k$.
    
    {\bf Claim.} There is a $Q>0$ such that for all $k$ $$\sum
    \limits_{U \in \, {\cal U}}(diam(U))^{d} > Q^{-1}\sum \limits_{U
    \in \, {\cal U}_{k}}(diam(U))^{d}.$$ Given the claim we see that
    the right hand side is the sum over all level-k clones.  We have
    seen above that as $k \to \infty$ this sum converges to $K'>0$. 
    Thus for $k$ large we have $$\sum \limits_{U \in \, {\cal
    U}_{k}}(diam(U))^{d} > \frac{K'}{2}.$$ Thus we have shown that
    $\eta' = Q^{-1}(K'/2)$ satisfies (\ref{eta}).  It remains to prove the claim.
    
    {\em Proof of claim.} It suffices to show this when ${\cal U}$
    contains only one clone $J$.  By lemma \ref{B}, there is a $Q>0$
    such that $\| M^{k} \vec{v} \|_{1} \leq Q \|\vec{v} \|_{1}$ for
    all positive vectors $\vec{v}$ in $\mathbb{R}^{n}$.  Let $J$ be a
    clone of type $i \in \, \{1, \ldots, n\}.$ Then $\vec{v} =
    (diam(J))^{d}\vec{e_{i}}$ is the $d$-quantity of the clone.  The
    clone $J$ subdivided $k$ times forms the collection of clones
    ${\cal J}^{(k)}$ with $d$-quantity $\vec{v}\left({\cal
    J}^{(k)}\right) = M^{k}\vec{v}\left( {\cal J} \right)$.  Then by
    lemma \ref{B} we have $$\sum \limits_{U \in \, {\cal
    J}^{(k)}}(diam(U))^{d} = \sum \limits_{i=1}^{n}
    |(M^{k}\vec{v})_{i}| = \|M^{k}\vec{v}\|_{1} \leq Q \| \vec{v}
    \|_{1} = Q(diam(J))^{d}.$$ This proves the claim.
    
    This proves that the Hausdorff measure of $C$ is finite and
    non-zero in this dimension.  Therefore, by lemma \ref{A}, we see
    that $d$ is the Hausdorff dimension.
    
    Let $\vec{v}$ be the vector where $v_{i}$ is the Hausdorff measure
    of the level-0 model $A_{i}$.  The $i,j$ entry of the transposed
    matrix $(M_{d})^{T}$ is the sum of the $d$-powers of the inverse
    scale factors of the clone-model maps from the type $j$ clones to
    the type $i$ models.  Then the $i$ entry of $(M_{d})^{T}\vec{v}$
    is $$\sum \limits_{j=1}^{n}(M_{d})^{T}_{ij}v_{j}.$$ By lemma
    \ref{multiply}, this is the sum of the Hausdorff $\delta$-measures
    of all the level-1 clones contained in model $i$.  Hausdorff
    measure is $\sigma$-additive so the Hausdorff $\delta$-measure of
    model $A_{i}$ must equal the measure of the level-1 clones
    contained in $A_{i}$.  Therefore $$(M_{d})^{T}\vec{v} = \vec{v}.$$
    This implies $\vec{v}$ is an eigenvector with eigenvalue 1.
    \end{proof}

\section {Measure Linearity on Multi-Model Cantor Sets}
\label{Linearity}

In this section we shall generalize some results by Cooper \cite{C} on
clone Cantor sets to the case of multi-model Cantor sets.  The main
result is theorem \ref{10.6} which states that given a bilipschitz map 
between multi-model Cantor sets then the map is measure linear on some 
clone.\\


\noindent {\bf Definition.}  Given a subset $A$ of a
clone Cantor set $C$, the {\em separation of $A$} is the minimum
distance between $A$ and $C-A$.
$$sep(A) = inf \{ d_{C}(x,y) : x \in \,  A,\, y \in \,  C-A \}.$$
By convention, $sep(C) = \infty$.  The {\em relative separation of
$A$} is the quotient of the separation of $A$ by the diameter of $A$: 
$$rel(A) = \frac{sep(A)}{diam(A)}.$$
    
As we discuss concepts related to Cantor sets we will often want to
consider clones after ``rescaling''.  Applying a similarity map with
stretch factor $K$ to clones increases the diameter of the clones by a
factor $K$ and also increases the diameter of the gaps between clones
by the same factor.  Thus the quotient of these values remains
constant.  Relative separation is therefore an invariant factor under
rescaling.\\


\noindent {\bf Definition.} The {\em minimum separation of C} is
denoted $$\alpha(C)=min\{sep(A) : A \mbox{ is a level-1 clone of
}C\}.$$

We now generalize some of the results from a paper written by Cooper
\cite{C} on clone Cantor sets.  We have made the necessary changes and
additions to make these work in the case of the multi-model Cantor
set.  We have rewritten the statements of the lemmas and the main
theorem from this paper and have summarized each.  We have outlined
the main points of the proofs.  The reader should refer to \cite{C} 
for further details.

Lemma \ref{10.2} includes a change of wording to signify
that $\tau$ is a map from a clone to the level-0 model of the same
type.  Here we use the definition of clone-model map given in section
\ref{multi-model}.  Lemmas \ref{10.3}, \ref{10.4} and \ref{10.5} now
involve an extra case where a clopen is not contained in any model and
theorem \ref{10.6} follows from these.


We refer to two quantities, dependent on certain
parameters, as being {\em approximately equal} if their ratio is bounded
above and below by positive constants independent of those 
parameters.

\begin{lemma}[Lemma 10.1 from \cite{C}]
    \label{10.1}
    Every clopen in a multi-model Cantor set is a finite union of clones.
\end{lemma}

\begin{lemma}[Lemma 10.2 from \cite{C}]
    \label{10.2}
    Given a multi-model Cantor set $C$ there is a constant $\xi(C) >
    1$ such that for every clone $A$ of $C$ with level of $A \geq 1$
    we have $$\frac{1} {\xi(C)} \leq rel(A) \leq \xi(C).$$
\end{lemma}

This lemma tells us that the separation of a clone approximates its
diameter.  The replacement of the clone map by a clone-model map
suffices to make the proof from \cite{C} apply to the multi-model
case.  We now sketch this.

There are only a finite number of clones contained in the Cantor set
of more than a certain diameter.  Two points are chosen, of minimal
distance apart, one inside clone $A$ and the other in $C-A$. 
Rescaling by a clone-model map does not change relative separation, so
we can assume the two points are contained in different level-1 clones
and that the distance between the two points is no less than the
minimum separation of $C$.  This now gives a lower bound on the
diameter of $D$, the smallest clone which properly contains $A$. 
There are only a finite number of possibilities for $D$ due to the
lower bound on its diameter.  Therefore there are a finite number of
options for $A$, as $A$ is at relative level 1 to $D$.  Since there
are only a finite number of choices for $A$ then there are only a
finite number of possibilities for $rel(A)$. \hfill $\square$

\begin{lemma}[Lemma 10.3 from \cite{C}]
    \label{10.3}
    Given a multi-model Cantor set C there is a constant $\beta$ with
    the following property.  Let B be any subset of C of positive
    diameter and let D be the smallest clone of C containing B (or
    choose D=C if there is no such clone).  Then $$diam(B) \leq
    diam(D) \leq
    \beta diam(B).$$
\end{lemma}

This lemma says that the diameter of any subset is approximately
the same as the diameter of the smallest clone which contains it.

Given a multi-model Cantor set $C$, we note that not every subset of
$C$ is contained in a clone.  However, if $B$ is
not contained in any clone then $diam(B)\geq \alpha(C)$.  Taking $D$
to be the whole Cantor set then we already have $$\frac{diam(D)}
{diam(B)} \leq \frac{diam(C)} {\alpha(C)}$$ and the result follows as
below.

If $B$ is contained in a clone then rescale $D$ onto the appropriate
level-0 model $A_{j}$, using a clone-model map $\tau$.  After
rescaling, $B$ must contain points in at least two level-1 clones by
choice of $D$.  Therefore the diameter of $B$ must be at least as
great as the minimum separation of $C$.  Now rescaling multiplies the
diameter of all clones by the same constant so the ratio of the diameter of $D$
to $B$ is the same as the ratio of the diameters after rescaling. 
This is then no greater than the ratio of the diameter of $A_{j}$ to
the minimum separation of $C$.  Therefore $$\frac{diam(D)}{diam(B)} =
\frac{diam(\tau(D))} {diam(\tau(B))} \leq \frac{diam(A_{j})}
{\alpha(C)} \leq \frac{diam(C)}{\alpha(C)}.$$ So if we let $\beta =
diam(C)(\alpha(C))^{-1}$ we get the desired result. \hfill $\square$

\begin{lemma}[Lemma 10.4 in \cite{C}]
    \label{10.4} Given $K > 1$, $\epsilon > 0$ and multi-model Cantor sets
    C, $C'$ there is a constant M with the following property.  Let f
    be any K-bilipschitz map of C onto a clopen in $C'$, and suppose
    that $rel(f(C)) \geq \epsilon$.  Let A be any clone of C. If $P'$
    is the smallest clone of $C'$ containing f(A) and $\tau'$ is the
    clone map taking $P'$ onto $C'$ then there are at most $M=M(C, C',
    K, \epsilon)$ possibilities for the image $\tau'(f(A))$ as a
    subset of $C'$.  This is also true if $f(A)$ is not
    contained in any clone of $C'$ with $\tau'$ taken to be the 
    identity map.

\end{lemma}

This lemma says that given a $K$-bilipschitz map from $C$ to $C'$ and a
clone $A$ in $C$, up to rescaling in $C'$, there are only
a finite number of possibilities for the image of $A$.

The proof given in \cite{C} does not consider the case where the
clopen $f(A)$ is not contained in any clone.  Suppose the clopen
$f(A)$ is not contained in any clone.  Then there exist points $x$ and
$y$ in $f(A)$ such that $x$ and $y$ are in different models, i.e.
level-0 clones.  This implies $diam(f(A)) \geq \alpha(C')$.  The rest
of the proof in \cite{C} then applies but with $\tau'$ replaced by the
identity map.  We now sketch this proof.

The relative separation is not affected by rescaling.  Given a clone
$A$ in $C$, then the diameter of $f(A)$ is bounded above by a factor
of $K$ times the diameter of $A$ because $f$ is $K$-bilipschitz.  The
separation of $f(A)$ is bounded below by a factor of $K^{-1}$ times
the separation of $A$ for the same reason.  Hence the relative
separation of $f(A)$ is bounded below by some positive constant.  We
then rescale $f(A)$ using a clone-model map $\tau'$ and find that the
diameter of $\tau'f(A)$ is bounded below by the minimum separation of
models in $C'$ and as relative separation is unchanged by rescaling we
find the separation of $\tau'f(A)$ is also bounded below.  We then
choose $B$ to be the biggest clone in $\tau'f(A)$ and let $D$ be the
smallest clone which properly contains $B$ (or all of $C'$).  The
diameter of $D$ is bounded below by the separation of $\tau'f(A)$ and
so there are only finitely many posibilities for $D$.  Clone $B$ is
one level higher than clone $D$ and hence there are only finitely many
possibilities for $B$.  The number of maximal clones in $\tau'f(A)$ is
bounded above because the diameter of each is bounded below.  Since 
there are a finite number of possibilities for the image of each of 
these clones we have a finite number of possibilities for
$\tau'f(A)$. \hfill $\square$

\begin{lemma}[Corollary 10.5 in \cite{C}]
    \label{10.5}
    Given $K>1$, $\epsilon > 0$ and multi-model Cantor sets C, $C'$ there is a 
    finite set S of positive real numbers with the following property.  
    Suppose that f is any K-bilipschitz map of C onto a clopen in $C'$, 
    and suppose that $rel(f(C)) \geq \epsilon$.  If B is any clone of level 
    n in $C$ containing a clone D of level n+1 then 
    $$\frac{{\cal H}_{\delta}(f(D))}{{\cal H}_{\delta}(f(B))}\, \in 
    \, S.$$
\end{lemma}

This may be interpreted as saying that non-linearity, at the level of measure
theory, comes in discrete ``packets''.  A bilipschitz map need not
distort the measure of all clones by the same multiplicative factor 
but, if $D$ is a clone contained in $B$ of relative level 1, then the ratio 
of the masses of the images of $B$ and $D$ can be one of only finitely 
many possible numbers.

In generalizing the proof to the case of multi-model Cantor sets we
need to consider the possibility that $f(B)$ is not contained in any
clone.  By lemma \ref{10.4} we already know there are a finite number
of possibilities for $f(B)$.  We need to show there are only finitely
many possibilities for $f(D)$ and the result will follow.

We consider two cases below.  If $f(D)$ is contained in a clone then
case 1 applies.  If $f(D)$ is not contained in any clone then
case 2 applies.

{\bf Case 1.} Suppose $f(D)$ is contained in some clone $Q'$.  We know
that $$\frac{diam(D)}{diam(B)} = \frac{diam(E_{i})}{diam(A_{j})},$$
where $A_{j}$ is a model and $E_{i}$ is a level-1 clone.  Hence this
quotient is equal to one of a finite number of possible positive scale
factors and so is bounded away from 0.  Now $f(B)$ is not contained in
a model so $diam(f(B)) \geq \alpha(C')$ and hence $diam(B) \geq
K^{-1}\alpha(C')$ as $f$ is $K$-bilipschitz.  We have that $diam(Q') $
approximates $diam(f(D))$ by lemma \ref{10.3} and so approximates
$diam(D)$ because $f$ is $K$-bilipschitz.  Hence $$ \frac{diam (Q')}{
K^{-1}\alpha(C')} \, \geq \, \mbox{ approximately} \,\,
\frac{diam(D)}{diam(B)}$$ and so is bounded away from 0 and thus
$diam(Q')$ is bounded away from 0.  There are only a finite number of
clones this large so we have only finitely many choices for $Q'$.  We
map $Q'$ onto the appropriate level-0 model $A_{j}$ by a clone-model
map $\phi$, (using $Q'=A_{j}$ and $\phi =$ identity if level $Q' = 0$). 
Then $\phi(f(D))$ is one of only finitely many possibilities by lemma
\ref{10.4} and since there are only a finite number of choices for
$\phi$ then there are only finitely many possibilities for the clopen
$f(D)$.

{\bf Case 2.}  Suppose $f(D)$ is not contained in any clone.  Then by 
lemma \ref{10.4} we have only finitely many possibilities for 
$f(D)$.

If $f(B)$ is contained in a clone then the rest of the proof for lemma
\ref{10.5} follows from the observation that there are a finite number
of possibilities for the rescaled image of clone $B$, and the rescaled
image of clone $D$, where $D$ contains $B$ at relative level 1. 
Therefore there are only a finite number of values for $$\frac{{\cal
H}_{\delta} (\tau'f(D))} {{\cal H}_{\delta}(\tau'f(B))}.$$ We take
$\tau'$ to be the identity map if $f(B)$ is not contained in any
clone.

In either case the map $\tau'$ is a similarity map with stretch
factor $L$ we have, using lemma \ref{multiply}, $$\frac{{\cal
H}_{\delta}(\tau'f(D))}{{\cal H}_{\delta}(\tau'f(B))}=
\frac{L^{\delta}{\cal H}_{\delta}(f(D))}{L^{\delta}{\cal
H}_{\delta}(f(B))} = \frac{{\cal H}_{\delta}(f(D))}{{\cal
H}_{\delta}(f(B))}.$$ Thus there are only a finite number of
possibilities for this ratio. \hfill $\square$\\

\noindent {\bf Definition.} Given a $K$-bilipschitz map $f:C \to C'$,
the {\em mass ratio} ($MR$) of a clone $A$ in $C$ is defined to be $$MR(A) =
    \frac{{\cal H}_{\delta}f(A)}{{\cal H}_{\delta}(A)}.$$

Thus if $f$ is the identity then $MR(A)=1$.
The larger this number the ``greedier'' we consider the clone $A$ to 
be.

\begin{cor}[Assertion in \cite{C}]
    \label{massratio}
    With the hypothesis of lemma \ref{10.5} there is a finite set S of
    real numbers such that for every pair of clones $E \subset D$ in C
    with $level(E) = level(D) +1$ then $$\frac{MR(E)}{MR(D)} \, \in 
    \, S.$$
\end{cor}

Here we see that when a clone $D$ and a clone $E$ contained within it
at relative level 1 are mapped into a multi-model Cantor set with a
$K$-bilipschitz map then the quotient of the mass ratios of these clones
is restricted to some finite set.\\

    \noindent {\bf Definition.} A map $f(X,\mu_{X}) \to (Y,\mu_{Y})$
    between metric spaces is {\em measure linear}, i.e. has constant
    Radon-Nikodym derivative, if for some constant $K$,
    $$\mu_{X}(f(A)) = K\mu_{Y}(A)$$ for all measurable subsets $A$ 
    of $X$.
    
This is a weak generalization of linearity.  Linearity refers to 
metric whereas this concept refers to measure.  For Lebesgue measure 
on the real line, linear and measure linear are the same 
for continuous functions.

\begin{theorem}[Theorem 10.6 in \cite{C}]
    \label{10.6}
    Suppose that C, $C'$ are multi-model Cantor sets and that f is a
    $K$-bilipschitz map of C onto a clopen subset of $C'$.  Then $C$ 
    and $C'$ have the same Hausdorff dimension.  Furthermore there is
    a clopen A in C such that the restriction $f|_{A}$ of f to A is
    measure linear with respect to Hausdorff measure.
\end{theorem}

In other words, any $K$-bilipschitz map from one multi-model Cantor set to
another is measure linear on some clopen.

We sketch the proof from \cite{C}.  The map $f$ is $K$-bilipschitz so
the mass ratio function is bounded above.  Thus there exists a
supremum for the possible mass ratios.  Let $A$ be a clone contained
in $C$ with mass ratio within $\epsilon$ of this supremum.  First note
that if the mass ratios of all clones contained in $A$ were equal,
then we would already have a measure linear map on $A$.  Also if each
clone contained in $A$ had mass ratio smaller than that of $A$ then
the union of these clones would also have smaller mass ratio which is
a contradiction.  So assume there is a clone $E$ of minimum level,
with mass ratio greater than that of $A$.  We call this a ``greedy''
clone.  We then, by corollary \ref{massratio}, have only a finite
number of possibilities for the quotient of mass ratios of $E$ and the
smallest clone $D$ which contains it.

If the mass ratio of $D$ were larger than the mass ratio of $A$ we
would have chosen this for $E$ as it has lower level.  On the other
hand the mass ratio of $D$ cannot be smaller than the mass ratio of
$A$ else there exists a greedy clone at the same level as $D$. 
Therefore the mass ratio of $D$ must be equal to the mass ratio of
$A$.  This implies that the quotient of mass ratios of $E$ and $D$
must be greater than 1.

Thus $MR(E)/MR(A)$ is in a finite set of numbers bigger than 1.  This
implies that the mass ratio of $E$ is larger than the supremum of mass
ratios, which is of course a contradiction and so proves that there
could not have been a ``greedy'' clone.  Therefore the map must have
been measure linear on $A$. \hfill $\square$

In \cite{C}, Cooper defines various invariants of clone Cantor
sets under bilipschitz maps.  We say two sets of real numbers $A$ and
$B$ are {\em similar} if there exist non-zero scalars $\alpha, \beta$
such that $$\alpha A \subset B \quad \mbox{and} \quad \beta B \subset
A.$$
Let $C$ be a multi-model Cantor set with Hausdorff measure ${\cal
H}_{\delta}(C)$ where $\delta$ is the Hausdorff dimension.  We define
the {\em clopen invariant} of an arbitrary subset $B$ of $C$ as the
similarity class of the countable set of real numbers $\{ {\cal
H}_{\delta}(A): A \mbox{ a clopen in } B \}$.

\begin{lemma}
    \label{C}
    Given a multi-model Cantor set with irreducible matrix and
    any two clones, $A$ and $B$, then the clopen invariants of $A$ and
    $B$ are equal.
\end{lemma}

\begin{proof}
    The clone-model map from A to the model $A_{j}$ shows that $A$ has
    the same clopen invariant as model $A_{j}$.  Similarly $B$ has the
    same clopen invariant as $A_{i}$ for some $i$.  A scaled down copy
    of $A_{j}$ is contained in $A_{i}$ for all i hence the clopen
    invariant for the model $A_{j}$ is similar to that for every other
    model.
\end{proof}

We can now define the clopen invariant of a multi-model Cantor
set.  Let $C$ be a multi-model Cantor set with Hausdorff measure ${\cal
H}_{\delta}(C)$ where $\delta$ is the Hausdorff dimension.  We define 
the {\em clopen invariant} of $C$ as the similarity class of the 
countable set of real numbers $\{ {\cal H}_{\delta}(A): A \mbox{ a 
clopen in } A_{j} \}$ for some $j$.  The choice of $j$ is irrelevant 
by lemma \ref{C}.

\begin{cor}[Corollary 10.7 from \cite{C}]
    \label{10.7}
    If $C$, $C'$ are multi-model Cantor sets and there is a 
    bilipschitz map of $C$ onto a clopen in $C'$, then $C$ and $C'$ 
    have the same clopen invariant.
\end{cor}

    By theorem \ref{10.6}, there is a clopen, hence a clone $A$ in $C$
    with a measure linear map of $A$ onto a clopen in $C'$.  The
    clone-model map from $A$ to $A_{j}$ for some $j$ shows that $A$
    has the same clopen invariant as $A_{j}$ and hence $C$ by lemma
    \ref{C}.  Let $\alpha$ be the Radon-Nikodym derivative of the
    measure linear map of $A$ into $C'$.  The image of a clopen $B$ in
    $A$ is a clopen in $C'$ with $\alpha$ times the measure of $B$. 
    By the same arguement as above, this clopen invariant is the same
    as the clopen invariant for $C'$.  Using that the inverse of a
    bilipschitz map is bilipschitz, we obtain the inverse
    relationship.

\noindent Primary subject: 37D40 Dynamical systems of geometric origin.\\
Secondary subject: 28A80 Fractals\\

\noindent E. Cockerill\\
Department of Mathematics\\
University of California\\
Santa Barbara, CA 93106\\
lilith@math.ucsb.edu

\end{document}